\documentclass[a4paper, 12pt, leqno]{amsart}

\usepackage{amsmath,amssymb,amsthm}
\usepackage{graphicx}
\usepackage{cite,mathrsfs}

\newtheorem{theorem}{Theorem}[section]
\newtheorem{lemma}[theorem]{Lemma}

\theoremstyle{definition}

\theoremstyle{remark}
\newtheorem{remark}[theorem]{Remark}

\numberwithin{equation}{section}

\begin{document}

\title{On the theory of $M$-functions}

\author{Kohji Matsumoto}
\address{K. Matsumoto: Graduate School of Mathematics, Nagoya University, Chikusa-\
ku, Nagoya 464-8602, Japan}
\email{kohjimat@math.nagoya-u.ac.jp}


\maketitle
\baselineskip 16pt

\begin{center}
{\it Dedicated to Professor Yasutaka Ihara \\on the occasion of his eightieth birthday}
\end{center}

\begin{abstract}
We survey the value-distribution theory of zeta and $L$-functions, originated by
H. Bohr, and developed further by Y. Ihara and others recently.
\end{abstract}
\bigskip

Our aim is to discuss the theory of {\bf $M$-functions}, which are certain kind of
``density functions" which describes the value-distribution of $L$-functions.

\section{The theory of Bohr and Jessen}

The name ``$M$-functions" is due to Ihara \cite{Iha08}.
His original motivation is to study the value-distribution of $(L'/L)(s,\chi)$, the
logarithmic derivative of Dirichlet or Hecke $L$-functions (on number fields, and on function 
fields).
This motivation is a natural extension of his study on Euler-Kronecker constants
(see \cite{Iha06}, \cite{Iha10}, \cite{IMS09}).

However, the primitive form of $M$-functions already appeared in the work of
Bohr and his colleagues in 1930s, on the value-distribution of the Riemann
zeta-function $\zeta(s)$.

First we recall the theory of Bohr and Jessen \cite{BJ3032} on the value-distribution of the Riemann 
zeta-function $\zeta(s)$.

$R(\subset\mathbb{C})$ : a rectangle with the edges parallel to the axes,

$\sigma > 1/2$,

$\mu_d$ : the usual $d$-dimensional Lebesgue measure,

$V_{\sigma}(T,R,\zeta)=\mu_1(\{t\in[-T,T]\;|\;\log\zeta(\sigma+it)\in R\})$.

\noindent
Then we have

\begin{theorem}\label{thm-BJ}
{\rm (the Bohr-Jessen limit theorem)}

{\rm(I)} When $T\to\infty$, the limit
$$
W_{\sigma}(R,\zeta)=\lim_{T\to\infty}\frac{1}{2T}V_{\sigma}(T,R,\zeta)
$$
exists.

{\rm(II)} There exists a function $\mathcal{F}_{\sigma}(\cdot,\zeta)$ defined 
on $\mathbb{C}$,
which is continuous and everywhere non-negative, for which
$$
W_{\sigma}(R,\zeta)=\lim_{T\to\infty}\frac{1}{2T}V_{\sigma}(T,R,\zeta)=
\int_R \mathcal{F}_{\sigma}(w,\zeta)|dw|
$$
{\rm(}where $|dw|=(2\pi)^{-1}dudv$ for $w=u+iv${\rm)} holds.
\end{theorem}

We may say that $W_{\sigma}(R,\zeta)$ is the {\bf probability} of how many values of 
$\log\zeta(s)$ on the vertical line $\Re s=\sigma$ belong to $R$, and 
$\mathcal{F}_{\sigma}(\cdot,\zeta)$ is the {\bf density function} of this probability.

To understand the behavior of $W_{\sigma}(R,\zeta)$, it is desirable to obtain some
explicit construction of $\mathcal{F}_{\sigma}(\cdot,\zeta)$. 
Bohr and Jessen themselves gave a construction, but
in a rather complicated way.

The way of the proof of Theorem \ref{thm-BJ} is as follows.
First, the Euler product gives:
$$\log\zeta(\sigma+it)=-\sum_{n=1}^{\infty} \log(1-p_n^{-\sigma-it}),$$
where $p_n$ denotes the $n$th prime number.
Consider the finite truncation
$$f_N(\sigma+it)=-\sum_{n=1}^{N} \log(1-p_n^{-\sigma-it})
=-\sum_{n=1}^{N} \log(1-p_n^{-\sigma}e^{-it\log p_n}).
$$
To analyze the properties of $f_N(\sigma+it)$, we introduce the associated mapping
$S_N:[0,1)^N \to \mathbb{C}$ defined by
$$
S_N(\theta_1,\ldots,\theta_N)=-\sum_{n=1}^{N} \log(1-p_n^{-\sigma}e^{2\pi i\theta_n})
\quad(0\leq\theta_n<1).
$$

The proof of Bohr and Jessen is based on:

$\bullet$ A mean value theorem for $\zeta(s)$ (necessary to show that
$f_N(\sigma+it)$ approximates $\log\zeta(\sigma+it)$ in a certain mean value sense),

$\bullet$ The fact that $\{\log p_n\}_{n=1}^{\infty}$ is linearly independent over $\mathbb{Q}$,

$\bullet$ The geometry of the auxiliary mapping $S_N$, 
that is, 
each term $-\log(1-p_n^{-\sigma}e^{2\pi i\theta_n})$ describes a convex curve when 
$\theta_n$ moves from 0 to 1.   

Bohr and Jessen \cite{BJ29} developed a detailed theory on the
``sums'' of convex curves, which is essentially used in their proof of
Theorem \ref{thm-BJ}.

Later Jessen and Wintner \cite{JW35} published an alternative proof of the Bohr-Jessen 
theorem, which is more Fourier theoretic.
In their proof, instead of convexity, a certain inequality is used.
This inequality, now called the Jessen-Wintner
inequality, is also 
connected with geometric properties of curves.    We will explain later what is the
Jessen-Wintner inequality.

\begin{remark}
The Bohr-Jessen theorem is for $\log\zeta(s)$.    An analogous result
for $(\zeta'/\zeta)(s)$ was shown by Kershner and Wintner \cite{KW37}.
\end{remark}

If one tries to generalize Theorem \ref{thm-BJ} to more general zeta or $L$-functions, 
we encounter the difficulty that
the corresponding geometry becomes more complicated
(for example, the convexity is {\bf not valid} in general). 

But still,
Part (I) of the Bohr-Jessen theorem (the existence of the limit
$W_{\sigma}(R,\zeta)=\lim_{T\to\infty}(2T)^{-1}V_{\sigma}(T,R,\zeta)$) has been extended 
to the
case of quite general zeta-functions which have Euler products (\cite{Mat90} \cite{Mat92}),
by invoking

$\bullet$ Prokhorov's theorem (in \cite{Mat90}), or

$\bullet$ L{\'e}vy's convergence theorem (in \cite{Mat92}).

The proofs are very analytic (or better to say, probabilistic) and do not use geometric properties (such as the convexity).

However, to prove Part (II), that is the existence of the density function
$\mathcal{F}_{\sigma}(\cdot,\zeta)$ satisfying
$$
W_{\sigma}(R,\zeta)=\lim_{T\to\infty}\frac{1}{2T}V_{\sigma}(T,R,\zeta)=
\int_R \mathcal{F}_{\sigma}(w,\zeta)|dw|,
$$
it seems that the convexity, or the Jessen-Wintner inequality, is essentially necessary.
Therefore the analogue of the above formula was formerly proved only for
the case when the attached curve is convex:

$\bullet$ Dirichlet $L$-functions (Joyner \cite{Joy86}),

$\bullet$ Dedekind zeta-functions attached to Galois number fields (\cite{Mat92}).

But there are some recent developments in the non-convex case, which we will report 
later.

\section{Ihara's work and related results}

So far we discussed the value-distribution of 
zeta or $L$-functions when $t=\Im s$ varies.
But it is also possible to study the value-distriution from some different point
of view.
For example, for Dirichlet or Hecke $L$-functions $L(s,\chi)$, we may consider the 
{\bf modulus aspect}.
\bigskip

Ihara \cite{Iha08} studied the behavior of $(L'/L)(s,\chi)$ from this aspect, and proved
the limit formula
$$
{\rm Avg}_{\chi}\Phi\left(\frac{L'}{L}(s,\chi)\right)=\int_{\mathbb{C}}
M_{\sigma}(w)\Phi(w)|dw|
$$
for a certain average with respect to $\chi$, where $\Phi$ is a test function and
$M_{\sigma}(\cdot)$, the {\bf $M$-function} for $L'/L$, is 
the (explicitly constructed) density function.

A typical meaning of Ihara's average ${\rm Avg}_{\chi}$ (in the rational number field
case) is:
$$
{\rm Avg}_{\chi}\phi(\chi)=\lim_{m\to\infty}\frac{1}{\pi(m)}\sum_{p\leq m}
\frac{1}{p-2}{\sum_{\chi({\rm mod}\; p)}}^{\!\!\!\!\!*} \phi(\chi),
$$
where $p$ runs over prime numbers, 
$\pi(m)$ is the number of prime numbers up to $m$,
and ${\sum_{\chi({\rm mod}\; p)}}^*$ means the sum
on primitive characters mod $p$.    Ihara's results in \cite{Iha08} are:

$\bullet$ In number field case, Ihara proved the formula (for any continuous function 
$\Phi$) in the region $\sigma>1$.    

$\bullet$
In the function field case, he proved the
same formula in wider region (such as $\sigma>1/2$, or $\sigma>3/4$)
for some special choices of $\Phi$, by using the ``proved'' Riemann Hypothesis.

As another average,
it is also possible to consider the character of the form
$\chi_{\tau}(p)=p^{-i\tau}$.     Then the associated $L$-function is
$$
\prod_p (1-\chi_{\tau}(p)p^{-s})^{-1}=\prod_p (1-p^{-s-i\tau})^{-1}
=\zeta(s+i\tau).
$$
The average associated with this type of character is
$$
{\rm Avg}_{\chi}\phi(\chi_{\tau})=\lim_{T\to\infty}\frac{1}{2T}\int_{-T}^T
\phi(\chi_{\tau})d\tau.
$$
For this average we can also prove \eqref{1} and \eqref{2}, which especially gives a
generalization of the result of Bohr and Jessen (their result is the case when $\Phi$
is the characteristic function of $R$).
In particular, $\mathcal{F}_{\sigma}(w,\zeta)$ of Bohr-Jessen is a special case of
$M$-functions.

Now let $L(s,\chi)$ be Dirichlet's, and recall Ihara's identity:
\begin{align}\label{1}
{\rm Avg}_{\chi}\Phi\left(\frac{L'}{L}(s,\chi)\right)=\int_{\mathbb{C}}
M_{\sigma}(w)\Phi(w)|dw|
\end{align}
and its ``log"-analogue:
\begin{align}\label{2}
{\rm Avg}_{\chi}\Phi\left(\log L(s,\chi)\right)=\int_{\mathbb{C}}
\mathcal{M}_{\sigma}(w)\Phi(w)|dw|.
\end{align}
Using certain mean value results, we can go into the region $1/2<\sigma\leq 1$.
\begin{theorem}\label{IM}
{\rm(Ihara and Matsumoto \cite{IM11} \cite{IM14})}
For $1/2<\sigma\leq 1$, both \eqref{1} and \eqref{2} hold with (explicitly constructed)
density functions $M_{\sigma}(w)$ and $\mathcal{M}_{\sigma}(w)$, for any $\Phi$ which is
(i) any bounded continuous function, or
(ii) the characteristic function of either a compact subset of $\mathbb{C}$ or the
complement of such a subset.
\end{theorem}

How to construct the $M$-function?
Here we explain the method in the $\log$ case.

Let $P$ be a finite set of primes, and let $N=|P|$.    Put
$$
L_P(s,\chi)=\prod_{p\in P}\left(1-\chi(p)p^{-s}\right)^{-1}.
$$
First we construct the density function $\mathcal{M}_{\sigma,P}(w)$ for which
$$
{\rm Avg}_{\chi}\Phi\left(\log L_P(s,\chi)\right)=\int_{\mathbb{C}}
\mathcal{M}_{\sigma,P}(w)\Phi(w)|dw|
$$
holds for any continuous $\Phi$.

Let $T=\{t\in\mathbb{C}\;|\;|t|=1\}$, and
define the auxiliary mapping 
$g_N:T^N \to\mathbb{C}$ by
$$
g_N((t_p)_{p\in P})=-\sum_{p\in P}\log(1-t_p p^{-\sigma}).
$$
(Note: This $g_N$ is essentially the same as the Bohr-Jessen auxiliary mapping
$S_N:[0,1)^N \to\mathbb{C}$.)

Using this $g_N$, and applying the orthogonality relation of characters, we find that
$$
{\rm Avg}_{\chi}\Phi\left(\log L_P(s,\chi)\right)=
\int_{T^N}\Phi(g_N((t_p)_{p\in P}))d^*T^N,
$$
where $d^*T^N$ is the normalized Haar measure on $T^N$.   Therefore our aim is to construct 
$\mathcal{M}_{\sigma,P}(w)$ for which
$$
\int_{T^N}\Phi(g_N((t_p)_{p\in P}))d^*T^N=\int_{\mathbb{C}}
\mathcal{M}_{\sigma,P}(w)\Phi(w)|dw|
$$
holds.

When $P=\{p\}$, we define
$$
\mathcal{M}_{\sigma,\{p\}}(w)=\frac{|1-r_p e^{i\theta_p}|^2}{r_p}\delta(r_p-p^{-\sigma}),
$$
where $r_p,\theta_p$ are determined by $w=-\log(1-r_pe^{i\theta_p})$ and 
$\delta(\cdot)$ denotes the Dirac delta distribution.

When $|P|\geq 2$ and $P=P'\cup\{p\}$, we define $\mathcal{M}_{\sigma,P}(w)$ 
recursively by the
convolution product
$$
\mathcal{M}_{\sigma,P}(w)=\int_{\mathbb{C}}\mathcal{M}_{\sigma,P'}(w')
\mathcal{M}_{\sigma,\{p\}}(w-w')|dw'|.
$$
This is a non-negative, compactly supported function, which satisfies the desired
property.

Next we have to show the existence of the limit
$$
\mathcal{M}_{\sigma}(w)=\lim_{|P|\to\infty}\mathcal{M}_{\sigma,P}(w).
$$
For this aim, we consider the Fourier transform
$$
\widetilde{\mathcal{M}}_{\sigma,P}(z)=\prod_{p\in P}\int_{\mathbb{C}}
\mathcal{M}_{\sigma,\{p\}}(w)\psi_z(w)|dw|,
$$
where $\psi_z(w)=\exp(i\Re(\overline{z}w))$.
Using the Jessen-Wintner inequality we can show that the right-hand side is
$O((1+|z|)^{-|P|/2})$, from which we can prove the existence of the limit
$$
\widetilde{\mathcal{M}}_{\sigma}(z)=
\lim_{|P|\to\infty}\widetilde{\mathcal{M}}_{\sigma,P}(z),
$$
and hence the existence of $\mathcal{M}_{\sigma}(w)$.

\begin{remark}
We can show the Dirichlet series expansion
$$
\widetilde{\mathcal{M}}_{\sigma}(z)=\sum_{n=1}^{\infty}\lambda_z(n)\lambda_{\overline{z}}(n)
n^{-2\sigma} \qquad (\sigma>1/2),
$$
where $\lambda_z(n)$ is defined by
$$
L(s,\chi)^{iz/2}=\sum_{n=1}^{\infty}\lambda_z(n)\chi(n)n^{-s}.
$$
Ihara \cite{Iha11} \cite{Iha12} studied a more general form
$$
\widetilde{\mathcal{M}}_{s}(z_1,z_2)=
\sum_{n=1}^{\infty}\lambda_{z_1}(n)\lambda_{z_2}(n)n^{-2s}
\qquad (\Re s>1/2),
$$ 
and proved various interesting properties.
\end{remark}

How general the test function $\Phi$ can be?     

In Theorem \ref{IM}, it is bounded
continuous, or the characteristic function of some compact subset, etc.

\begin{theorem}\label{IM2}
{\rm (Ihara and Matsumoto \cite{IM11b})}
If we assume the Generalized Riemann Hypothesis, the same type of limit theorem 
(for a little different definition of ${\rm Avg}_{\chi}$) holds
for any continuous $\Phi$ of at most exponential growth (that is, 
$\Phi(w)=O(e^{a|w|})$ with some $a>0$).
\end{theorem}

\begin{remark}
Therefore in the function field case, this theorem holds
unconditionally.    In this case a little weaker result was already obtained in
\cite{IM10}.
\end{remark}

\begin{remark}
To prove Theorem \ref{IM2}, the generalized form 
$\widetilde{\mathcal{M}}_{s}(z_1,z_2)$ (mentioned above) is necessary.
\end{remark}

\begin{remark}
An announcement of the above results of Ihara and the author appeared in \cite{IM09}.
\end{remark}

We list up some recent developments in the theory of $M$-functions.

$\bullet$
Mourtada and V. K. Murty \cite{MM15}
considered the average of
$(L'/L)(\sigma,\chi_D)$, where $\sigma>1/2$, $D$ is a fundamental discriminant and
$\chi_D$ is the associated real character, as $D\to\infty$,
and proved the same type of limit theorem.
Akbary and Hamieh \cite{AHarxiv} treated the cubic character case, and
Gao and Zhao \cite{GZarxiv} studied the quartic case.

$\bullet$
Suzuki \cite{Suz15} discovered that the $M$-function also
appears in the study of the vertical distribution of the zeros of certain functions
related with $\zeta(s)$.
Let $\xi(s)=\frac{1}{2}s(s-1)\pi^{-s/2}\Gamma(s/2)\zeta(s)$, $\omega>0$, and define
$$
A_{\omega}(s)=\frac{1}{2}(\xi(s+\omega)+\xi(s-\omega)),
B_{\omega}(s)=\frac{1}{2}i(\xi(s+\omega)-\xi(s-\omega)).
$$
Arrange the zeros of $A_{\omega}(s)$ (or $B_{\omega}(s)$) as $\rho_n=\beta_n+i\gamma_n$,
$\gamma_{n+1}\geq \gamma_n>0$.    Then it is known that the normalized imaginary part
$$
\gamma_n^{(1)}=\frac{\gamma_n}{2\pi}\log\frac{\gamma_n}{2\pi e}
$$
is well-spaced.    Suzuki doscovered that the "second order"
normalization
$$
\gamma_n^{(2)}=\left(\frac{\gamma_n}{2\pi}\log\frac{\gamma_n}{2\pi e}-n\right)
\rho_{\omega}^{-1/2}\frac{1}{2\pi}\log\frac{\gamma_n}{2\pi e},
$$
where $\rho_{\omega}=(2\pi^2)^{-1}\sum_{n=1}^{\infty}\lambda(n)^2n^{-1-2\omega}$,
is also well-distributed, and its law can be written by an integral involving
the $M$-function.

$\bullet$
Mine \cite{Min1} studied the $M$-function for Dedekind zeta-functions $\zeta_F(s)$
($F$: number field).
If $F$ is Galois, then, as mentioned before, the original argument of Bohr-Jessen can be
applied (\cite{Mat92}).
But in the non-Galois case, the situation is more difficult.
Mine noticed that the idea of Guo \cite{Guo96a} \cite{Guo96b} (for the distribution of
$(\zeta'/\zeta)(s)$) can be applied, and obtained the construction of $M$-function for
$(\zeta_F'/\zeta_F)(s)$, for any number field $F$.    Further generalization was done by
Mine \cite{Min2}.

$\bullet$
Recently Mine \cite{Min3} studied the $M$-function for zeta-functions of Hurwitz type
(that is, without Euler product).

\section{The value-distribution of automorphic $L$-functions}

How can we construct the theory of $M$-functions in the case of automorphic $L$-functions?
In this case, the attached curves are not always convex.

Let $f$ be a primitive form (i.e. normalized Hecke-eigen new cusp form)
of weight $k$ and level $N$, whose Fourier expansion is
$$
f(z)=\sum_{n\geq 1}\lambda_f(n)n^{(k-1)/2}e^{2\pi inz}.
$$
The corresponding $L$-function
$L_f(s)=\sum_{n\geq 1}\lambda_f(n)n^{-s}$
has the Euler product expansion
\begin{align*}
L_f(s)&=\prod_{p|N}(1-\lambda_f(p)p^{-s})^{-1}\prod_{p\nmid N}(1-\lambda_f(p)p^{-s}
+p^{-2s})^{-1}\\
&=\prod_{p|N}(1-\lambda_f(p)p^{-s})^{-1}\prod_{p\nmid N}
(1-\alpha_f(p)p^{-s})^{-1}(1-\beta_f(p)p^{-s})^{-1}.
\end{align*}

First consider the {\bf $t$-aspect}.  
Let
$$
V_{\sigma}(T,R,L_f)=\mu_1(\{t\in[-T,T]\;|\;\log L_f(\sigma+it)\in R\}).
$$

\begin{theorem}\label{t-aspect}
{\rm(Matsumoto and Umegaki \cite{MU2})}
For any $\sigma>1/2$, the limit
$$W_{\sigma}(R,L_f)=\lim_{T\to\infty}(2T)^{-1}V_{\sigma}(T,R,L_f)$$
exists, and it can be written as
$$
W_{\sigma}(R,L_f)=\int_R \mathcal{F}_{\sigma}(w,L_f)|dw|,
$$
where $\mathcal{F}_{\sigma}(\cdot,L_f)$ is a continuous, non-negative function (explicitly
constructed) on $\mathbb{C}$.
\end{theorem}

(A key of the proof) We have to show an
analogue of the Jessen-Wintner inequality for the automorphic case.

What is the {\bf original} Jessen-Wintner inequality?   Recall
$$
S_N(\theta_1,\ldots,\theta_N)=\sum_{1\leq n\leq N} z_n(\theta_n)\quad
\left(z_n(\theta_n)=
-\log(1-p_n^{-\sigma}e^{2\pi i\theta_n})\right).
$$
The original Jessen-Wintner inequality \cite{JW35} is the estimate
$$
\int_0^1 e^{i\langle w,z_n(\theta)\rangle}d\theta \ll p_n^{\sigma/2}|w|^{-1/2}
\quad (w\in\mathbb{C}),
$$
where
$\langle w,z_n(\theta)\rangle=\Re w \Re z_n(\theta)+\Im w \Im z_n(\theta)$.

In the automorphic case, instead of $z_n(\theta)$, we have to consider
$$
z_{f,p}(\theta)=-\log(1-\alpha_f(p)p^{-\sigma}e^{2\pi i\theta})
-\log(1-\beta_f(p)p^{-\sigma}e^{2\pi i\theta}).
$$
\begin{lemma}
[Inequality of Jessen-Wintner type, \cite{MU2}]
$$
\int_0^1 e^{i\langle w, z_{f,p}(\theta)\rangle}d\theta \ll_{\varepsilon} 
p^{\sigma/2}|w|^{-1/2}+p^{\sigma}|w|^{-1} \quad (w\in\mathbb{C})
$$
holds for any $p\in \mathbb{P}_f(\varepsilon)$, where
$$
\mathbb{P}_f(\varepsilon)=\{p:{\rm prime}\;|\;|\lambda_f(p)|>\sqrt{2}-\varepsilon\}.
$$
\end{lemma}

It is known that $\mathbb{P}_f(\varepsilon)$ is of positive density in the set of
all primes (M. R. Murty \cite{Mur83} in the full modular case, and in the book of
V. K. Murty and M. R. Murty \cite{MM97} in general case).    This is sufficient for our aim.

\begin{proof}[Proof of the lemma (sketch)]
Let $g(\theta)=\langle w, z_{f,p}(\theta)\rangle$.
Compute $g^{\prime}(\theta)$ and $g^{\prime\prime}(\theta)$.
Using the fact $|\lambda_f(p)|>\sqrt{2}-\varepsilon$, we can show that $[0,1)$ can
be divided into two subintervals $I_1$ and $I_2$, such that 
$|g^{\prime}(\theta)|$ is not small on $I_1$, while
$|g^{\prime\prime}(\theta)|$ is not small on $I_2$.
(That is, the geometric behavior of the curve $z_{f,p}(\theta)$ is ``not so bad".)
We apply the first derivative test on $I_1$, and the second derivative test on $I_2$.
\end{proof}

This lemma is the key of the proof of Theorem \ref{t-aspect}.
We omit how to deduce the theorem from the lemma.

We can also prove an analogous result for the $\gamma$-th symmetric power $L$-function 
for any $\gamma\in\mathbb{N}$, which is of the form
$$
L({\rm Sym}_f^{\gamma},s)=L_N({\rm Sym}_f^{\gamma},s)\prod_{p|N}({\rm certain
\;local \;factor \;at\;} p)
$$
where
$$
L_N({\rm Sym}_f^{\gamma},s)=\prod_{p\nmid N}\prod_{h=0}^{\gamma}(1-
\alpha_f^{\gamma-h}(p)\beta_f^h(p)p^{-s})^{-1}.
$$
Recently, connected with the Sato-Tate conjecture, there has been a big progress
on the study of symmetric power $L$-functions
(Barnet-Lamb et al. \cite{BGHT}).
Under certain plausible assumptions, using a result of \cite{BGHT} or its
quantitative version due to Thorner \cite{Tho14}, we can show an analogue of Theorem \ref{t-aspect} 
(that is, the existence of the associated $M$-function) for
$L({\rm Sym}_f^{\gamma},s)$
(see \cite{MU3}).

\bigskip

Secondly, the {\bf modulus aspect}.

Consider the twisted automorphic $L$-functions $L_f(s,\chi)$ whose local factor is defined by
$$
(1-\alpha_f(p)\chi(p)p^{-s})^{-1}(1-\beta_f(p)\chi(p)p^{-s})^{-1}.
$$
Lebacque and Zykin \cite{LZarxiv} obtained the formulas similar to
\eqref{1} and \eqref{2} for certain average of $L_f(s,\chi)$ with respect to
characters.

\bigskip

Thirdly, the {\bf level aspect}.    

So far there are two attempts:
Lebacque and Zykin \cite{LZarxiv}, and the author and Umegaki \cite{MU1}.

In the work of Bohr and Jessen, an essential fact is the linear independence of
$\{\log p_n\}$, and in the case of character-average we need the orthogonality
property of Dirichlet characters.

In the level-aspect case for automorphic $L$-functions, the corresponding tool is
Petersson's formula, which is used in both of the above articles.

Hereafter we explain the result of the author and Umegaki.

Consider the case $N=q^m$, where $q$ is a prime.

Define the $\gamma$-th (partial) symmetric power $L$-function by
$$
L_q({\rm Sym}_f^{\gamma},s)=\prod_{p\neq q}\prod_{h=0}^{\gamma}(1-
\alpha_f^{\gamma-h}(p)\beta_f^h(p)p^{-s})^{-1}.
$$
Assume:
(H1) $L_q({\rm Sym}_f^{\gamma},s)$ can be continued holomorphically to $\sigma>1/2$, and
$L_q({\rm Sym}_f^{\gamma},s)\ll N(|t|+2)$ in the strip $1/2<\sigma<2$;
(H2) There is no zero of $L_q({\rm Sym}_f^{\gamma},s)$ in the strip $1/2<\sigma\leq 1$.

Denote by $B_k(q^m)$ the set of all primitive forms of weight $k$ and level
$q^m$, consider certain weighted average on $B_k(q^m)$, and then take the limit
\begin{align*}
{\rm Avg}_{\rm prime}=\lim_{q\to\infty}
\;\;(m:{\rm fixed}),\;\;{\rm or}\;\;
{\rm Avg}_{\rm power}=\lim_{m\to\infty}
\;\;(q:{\rm fixed}).
\end{align*}

\begin{theorem}\label{symmetric}
{\rm (Matsumoto and Umegaki \cite{MU1})}
Let $2\leq k\leq 10$ or $k=14$, $\mu,\nu\in\mathbb{N}$ with $\mu-\nu=2$, and assume
{\rm (H1), (H2)} for the $\mu$-th and $\nu$-th symmetric power $L$-functions.   Then 
for any $\sigma>1/2$, there exists an explicitly constructed density function
$\mathcal{M}_{\sigma}:\mathbb{R}\to\mathbb{R}_{\geq 0}$ for which
\begin{align*}
&{\rm Avg}_{\rm prime}\Phi(\log L_q({\rm Sym}_f^{\mu},\sigma)
-\log L_q({\rm Sym}_f^{\nu},\sigma))\\
&= {\rm Avg}_{\rm power}\Phi(\log L_q({\rm Sym}_f^{\mu},\sigma)
-\log L_q({\rm Sym}_f^{\nu},\sigma))\\        
&=\int_{\mathbb{R}}\mathcal{M}_{\sigma}(u)\Phi(u)\frac{du}{\sqrt{2\pi}}
\end{align*}
holds for any $\Phi$ which is bounded continuous, or the characteristic function of 
either a compact subset of $\mathbb{R}$ or the
complement of such a subset.
\end{theorem}

\begin{remark}
We explain the reason why we study the "difference" of two $L$-functions.
It is possible to show that
\begin{align*}
\lefteqn{\log L({\rm Sym}_f^{\mu},\sigma)-\log L({\rm Sym}_f^{\nu},\sigma)}\\
&=-\sum_{p\neq q}\bigl(\log(1-\alpha_f^{\mu}(p)p^{-\sigma})+
\log(1-\beta_f^{\mu}(p)p^{-\sigma})\bigr).
\end{align*}
If we can take $\mu=1$, this is exactly $\log L_f(\sigma)$ (without the Euler
factor corresponding to $p=q$).    Therefore we could arrive at the theorem on
the value-distribution of $\log L_f(\sigma)$.

However, so far we cannot treat the case $\mu=1$.
To extend our result to the case $\mu=1$
is an important remaining problem.
\end{remark}

\end{document}